\RequirePackage{fix-cm}
\documentclass[smallextended]{svjour3}       
\smartqed  
\usepackage{graphicx}
\usepackage{graphics}
\usepackage{epstopdf}
\usepackage{amssymb,amsmath}
\usepackage{color}

\journalname{(under review in) Nonlinear Dynamics}
\begin{document}

\title{Novel Parameter Estimation Strategies for Time-Varying Systems via Real-Time Non-Linear Receding Horizon Control in Chaotic Environments
}

\titlerunning{Novel Parameter Estimation for Time-Varying Systems in Chaotic Environments}        

\author{Fei Sun         \and
        Kamran Turkoglu 
}


\institute{F. Sun \at
              Department of Aerospace Engineering, San Jos\'e State University, San Jose, CA 95192\\
              Tel.: +1-408-924-3965\\
              Fax: +1-408-924-3818\\
              \email{fei.sun@sjsu.edu}           
           \and
           K. Turkoglu \at
               Department of Aerospace Engineering, San Jos\'e State University, San Jose, CA 95192 \\
               \email{kamran.turkoglu@sjsu.edu}
}

\date{Received: date / Accepted: date}

\maketitle

\begin{abstract}
In this paper, based on real-time nonlinear receding horizon control methodology, a novel approach is developed for parameter estimation of time invariant and time varying nonlinear dynamical systems in chaotic environments. Here, the parameter estimation problem is converted into a family of finite horizon optimization control problems. The corresponding receding horizon control problem is then solved numerically, in real-time, without recourse to any iterative approximation methods by introducing the stabilized continuation method and backward sweep algorithm. The significance of this work lies in its real-time nature and its powerful results on nonlinear chaotic systems with time varying parameters. The effective nature of the proposed method is demonstrated on two chaotic systems, with time invariant and time varying parameters. At the end, robustness performance of the proposed algorithm against bounded noise is investigated.
\end{abstract}

\section{Introduction}
\label{intro}
Parameter estimation is a process of assessing unknown parameters with respect to a given limited amount of information. It is a procedure that is widely used in modeling and control of dynamical systems, where the applications range from biology to chemistry, physics and many others fields of science and engineering \cite{ho2010,wang2012,rao2012,lin2014,bellsky2014,ding2014,chi2014,wang2014nd}.

In real world applications, many systems do exhibit partially or completely unknown parameters. Knowledge about the time evolution of these parameters becomes the prerequisite to analyze, control, and predict the underlying dynamical behaviors. Thus, this topic has drawn great attention in various areas due to its theoretical and practical significance. For example, in the biological networks, it is important to estimate unknown protein-DNA interactions in the regulation of various cellular processes or detection of failures/anomalies. In aircraft/spacecraft dynamics, estimation of the unknown states determines the fine line between stability and instability.

Another significant area of interest for parameter estimation problems is chaotic systems. In many real-life problems, ranging from information sciences to life sciences, from systems biology to quantum physics \cite{pecora1990,elson1988,wu1993}, nonlinear systems exhibit the phenomenon of chaos. An important application of chaos control and synchronization is parameter estimation through adaptive control methodologies. In this case, the aim is to estimate the uncertainties as well as minimize the synchronization error \cite{ge2005,li2004,li2007,sun2009}. However, such adaptive control methodology is associated with the stability and the synchronization regime of systems. It is also a relatively conservative methodology that is constrained by several conditions such as persistent excitation or  linearly independence, to guarantee the convergence.

On the other hand, receding horizon control \cite{bryson1975,chen1982} is a branch of model predictive control methodology that aims to obtain an optimal feedback control law by minimizing the given performance index. The performance index of a receding-horizon control problem has a moving initial time and a moving terminal time, where the time interval of the performance index is finite. Since the time interval of the performance index is finite, the optimal feedback law can be determined even for a system that is not stabilizable. The receding horizon optimal control problem can deal with a broader class of control objectives than asymptotic stabilization \cite{ohtsuka1997}. The receding horizon control was originally applied to linear systems and then was extended to nonlinear systems \cite{thomas1975,kwon1983,mayne1990,mayne1990b}. Through its functionality, nonlinear receding horizon control (NRHC)\cite{mayne1990,mayne1990b} has made an important impact on industrial control applications and is being increasingly applied in process controls. Various advantages are known for NRHC, including the ability to handle time-varying and nonlinear systems, input/output constraints, associated plant uncertainties, and so on.

In recent years, many methods have been proposed for
parameter estimation in nonlinear systems. Some of them focused on using adaptive feedback
control algorithms to estimate unknown parameters of nonlinear systems. Huang \cite{huang2006pre} studied the adaptive synchronization with application to parameter estimation. Yu et al \cite{yu2007pre} proposed the linear independence conditions to ensure the parameter convergence based on the LaSalle invariance principle. However, most of these literatures \cite{huang2006pre,li2007,sun2009,yu2007pre,parlitz1996,sun2012} on adaptive estimation impose an assumption that the parameters to be estimated are constant or slowly time-varying. Moreover, the parameter estimation problem could be formed as an optimization problem, which leads to many intelligent optimization schemes:  Li et al. (\cite{li2006,li2010pre}) proposed the chaotic ant swarm algorithm and conducted parameter estimation tests by using the Lorenz system as an example. He et al. \cite{he2007} employed the particle swarm optimization method in parameter estimation. Lin et al. \cite{lin2014} proposed an oppositional seeker optimization algorithm with application to parameter estimation of chaotic systems.

Different from the above methods, in this paper, a method of parameter estimation for nonlinear systems is proposed based on real-time nonlinear receding horizon control (NRHC) methodology. With this approach, we provide a configuration which is especially applicable to chaotic and time varying systems. Here, the estimation procedure is reduced to a family of finite horizon optimization control problems. To avoid high computational complexity, the stabilized continuation method \cite{ohtsuka1998,ohtsuka1997} is employed, which is a non-iterative optimization procedure with moderate data storage capacity. Based on this method, the NRHC problem is then solved by the backward sweep algorithm \cite{bryson1975}, in real time. The algorithm itself is executable regardless of controllability or stabilizability of the system, which is one of the powerful aspects of the approach. Experimental results show that the real-time NRHC is applicable to the chaotic systems with unknown constant parameters as well as time-varying parameters. Furthermore, we explore the noise reduction of the proposed method by simulations.

In the light of these facts, the paper is organized as following: In Section-\ref{sec:prob_form}, the problem formulation, based on NRHC, is defined as an estimation routine. In Section-\ref{sec:sec: backwards_sweep}, brief background on previous work of Ohtsuka's \cite{ohtsuka1998,ohtsuka1997} is provided, and then stability analysis of this approach is discussed in Section-\ref{sec:sec:Conv_stab_analysis}. We demonstrate the power of NRHC as an estimation routine through specific applications on a chaotic system (in this case Lorenz oscillator \cite{lorenz1963}) with constant  (Section-\ref{sec:sec:app2chaotic_sys_const_params}) and time varying parameters (Section-\ref{sec:sec:app2chaotic_sys_TimeVar_params}). We also test and demonstrate robustness properties of NRHC algorithm in presence of noise, in Section-\ref{sec:sec:app2chaotic_sys_noise}. At the end, with the discussions and conclusions (Section-\ref{sec:conclusion}), we finalize the paper.

\section{Problem formulation}\label{sec:prob_form}
To demonstrate the parameter estimation routine of nonlinear systems, suppose that we are given the dynamical system representation as follows:
\begin{equation}\label{eq:GenNonLinSys}
\dot{x}=Ax+f(x)+D(x)\hat{\Theta}(t),\\
\end{equation}
where $x\in R^{n}$ is the state vector, $A\in R^{n\times n}$
and $f(x):R^{n}\rightarrow R^{n}$ are the linear coefficient matrix and nonlinear part of system presented in Eq.\eqref{eq:GenNonLinSys}, respectively. $D(x):R^{n}\rightarrow R^{n\times p}$ is a known function vector and $\hat{\Theta}\in R^{p}$
denotes the unknown parameters, where they can be constant or time-varying.

To formulate the parameter estimation problem, the system in Eq. (1) is considered as a drive/reference system. If we construct a driver-response configuration, the corresponding response system becomes

\begin{equation}
\dot{y}=Ay+f(y)+D(y)\Theta(t),\\
\end{equation}
where $y\in R^{n}$ is the state vector and $\Theta$ represents the estimated parameter. Here, functions $f(\cdot)$ and $D(\cdot)$ satisfy the global Lipschitz condition, therefore there exist
positive constants $\beta_1$ and $\beta_2$ such that
\begin{equation*}
\begin{split}
&\|f(y)-f(x)\|\le \beta_1\|y-x\|,\\
&\|D(y)-D(x)\|\le \beta_2\|y-x\|.
\end{split}
\end{equation*}
is satisfied.

In this specific formulation, the synchronization error $e(t) = y(t) -x(t)$ is defined to represent the difference between the drive system and the response system which is modeled as
\begin{equation}\label{eq:DrResp_sys_err}
\dot{e}(t)=Ae(t)+B(f(y(t))-f(x(t)))+D(y)\Theta-D(x)\hat{\Theta}.
\end{equation}
Also, the estimation error is denoted by $\bar{\Theta}(t)=\Theta(t)-\hat{\Theta}(t)$.

In order to utilize the real-time nonlinear receding horizon control method as an estimation routine, the following finite horizon cost function (performance index) is associated with the synchronization and estimation error:

\begin{equation}
\begin{split}
J=&\int_t^{t+T}L[y(\tau),x(\tau),\Theta(\tau),\hat{\Theta}(\tau)]{\rm d}\tau,\\
=&\int_t^{t+T}(e^TQe+\bar{\Theta}^TR\bar{\Theta}){\rm d}\tau.
\end{split}
\end{equation}
Here, $Q>0$ and $R>0$ are weighting matrices, affiliated with the state and estimation error, respectively. In this specific set-up, the performance index evaluates the performance from the present time-($t$) to the finite future-($t+T$), where $T$ is the terminal time or the horizon. The performance index is minimized for each time $t$ starting from $y(t)$. Thus, the present receding horizon control problem can be converted to a family of finite horizon optimal control problems on the $\tau$ axis that is parameterized by time $t$. The trajectory $y(t+\tau)$ starting from $y(t)$ is denoted as $y(\tau,t)$. Since the performance index of receding horizon control is evaluated over a finite horizon, the value of the performance index is finite even if the system is not stabilizable.

It is well known from literature that first order necessary conditions of optimality are obtained from the two-point boundary value problem (TPBVP) \cite{bryson1975} by computing the variations as the following

\begin{equation}\label{tpbvp}
\begin{split}
&y_\tau(\tau,t) = H_{\lambda}^T ,\quad y(0,t) = y(t), \\
&\lambda_\tau(\tau,t) = -H_y ^T,\quad \lambda(T,t) = 0, \\
&H_{\bar{\Theta}} = 0,
\end{split}
\end{equation}
In Eqs.\eqref{tpbvp}, $H$ is the Hamiltonian defined as
\begin{equation}
\begin{split}
H &= L + \lambda^{T}\dot{y}\\
&=(e^TQe+\bar{\Theta}^TR\bar{\Theta})+ \lambda^{T}[Ay+f(y)+D(y)\Theta].
\end{split}
\end{equation}
In this notation, $H_y$ denotes the partial derivative of $H$ with respect to $y$, and so on. According to this unique approach, the estimation error is calculated as
\begin{equation}\label{theta}
\bar{\Theta}(t) = \text{arg}\{H_{\bar{\Theta}}[y(t),\lambda(t),\bar{\Theta}(t),x(t)] = 0\}.
\end{equation}

%

In this context, the TPBVP is regarded as a nonlinear equation with respect to the costate at $\tau=0$ as
\begin{equation}
F(\lambda(t),y(t),T,t)=\lambda(T,t)=0.
\end{equation}
Since the nonlinear equation $F(\lambda(t),y(t),T,t)$ has to be satisfied at any time $t$, $\frac{{\rm d}F}{{\rm d}t}=0$ holds along the trajectory of the closed-loop system. The ordinary differential equation of $\lambda(t)$ can be solved numerically without applying any iterative optimization methods. However, numerical error in the solution may accumulate through the integration process in practice, and numerical stabilization techniques are required to correct the error. Therefore, the stabilized continuation method \cite{ohtsuka1998,ohtsuka1997} is introduced in this paper as follows:
\begin{equation}\label{cm}
\frac{{\rm d}F}{{\rm d}t}=-A_sF,
\end{equation}
where $A_s>0$ denotes any stable matrix and will enforce the exponential convergence of the solution. The horizon $T$ is defined as a smooth function of time $t$ such that $T(0)=0$ and $T(t)\rightarrow T_f$ as $t\rightarrow \infty$, where $T_f$ is the desired terminal time.


\subsection{Backward sweep algorithm:}\label{sec:sec: backwards_sweep}

In order to compute the estimation error, first, the differential equation of $\lambda(t)$ is integrated in real time. The partial differentiation of Eqs.\eqref{tpbvp} (with respect to time $t$ and $\tau$) converts the problem in hand into the following linear differential equation:
\begin{equation}\label{pd}
\frac {\partial}{\partial \tau}\begin{bmatrix}y_t-y_{\tau}\\
\lambda_t-\lambda_{\tau} \end{bmatrix}=\begin{bmatrix}G&-L\\
-K&-G^T \end{bmatrix}\begin{bmatrix}y_t-y_{\tau}\\
\lambda_t-\lambda_{\tau} \end{bmatrix}
\end{equation}
where $G=f_y-f_{\bar{\Theta}}H^{-1}_{\bar{\Theta}\bar{\Theta}}H_{\bar{\Theta}y}$, $L=f_{\bar{\Theta}}H^{-1}_{\bar{\Theta}\bar{\Theta}}f_{\bar{\Theta}}^T$, $K=H_{yy}-H_{y\bar{\Theta}}H^{-1}_{\bar{\Theta}\bar{\Theta}}H_{\bar{\Theta}y}$. Since the reference trajectory $x_t(t+\tau)=x_{\tau}(t+\tau)$, they are canceled in Eq.\eqref{pd} and data storage is reduced.

The derivative of the nonlinear function $F$ with respect to time is rewritten by
\begin{equation}\label{df}
\frac{{\rm d}F}{{\rm d}t}=\lambda_t(T,t)+\lambda_{\tau}(T,t)\frac{{\rm d}T}{{\rm d}t}.
\end{equation}

To reduce the computational cost, the backward-sweep algorithm is employed at this point and the relationship between the costate and other variables is expressed as:
\begin{equation}\label{relation}
\lambda_t-\lambda_\tau=S(\tau,t)(y_t-y_\tau)+c(\tau,t),
\end{equation}
where
\begin{equation}\label{sc}
\begin{split}
S_{\tau}&=-G^TS-SG+SLS-K, \\
c_{\tau}&=-(G^T-SL)c,
\end{split}
\end{equation}
In Eqs.\eqref{sc}, due to the terminal constraint on $\tau$-axis, the following conditions hold
\begin{equation}\label{sct}
\begin{split}
S(T,t)&=0,\\
c(T,t)&=H_y^T\mid_{\tau=T}(1+\frac{{\rm d}T}{{\rm d}t})-A_sF.
\end{split}
\end{equation}

Thus, the differential equation of $\lambda(t)$ is obtained in real time as follows:
\begin{align}\label{dl}
\frac{{\rm d}\lambda(t)}{{\rm d}t}=-H_y^T+c(0,t).
\end{align}

At each time $t$, the Euler-Lagrange equations Eqs.\eqref{tpbvp} are integrated forward along the $\tau$ axis. Eqs.\eqref{sc} are integrated backward with terminal conditions expressed in Eqs.\eqref{sct}. Then the differential equation of $\lambda(t)$ is integrated for one step along the $t$ axis so as to minimize the estimation error from Eq.\eqref{theta}. The estimated parameters are derived from the difference between the true values and the estimation errors. If the matrix $H_{\bar{\Theta}\bar{\Theta}}$ is nonsingular, the algorithm is executable regardless of controllability or stabilizability or the system.

\subsection{Stability analysis of NRHC estimation problem}\label{sec:sec:Conv_stab_analysis}

In this section, the stability of the closed-loop system by using the NRHC strategy is briefly analyzed. The candidate Lyapunov function is constructed in the form of
\begin{equation}
V=\frac{1}{2}e^Te.
\end{equation}
Here, clearly, $V(0)=0$ and $V>0$ for all $e\neq 0$, thus, $V$ is a Lyapunov function.

The time derivative of $V$ along the trajectory is obtained by
\begin{equation}
\begin{split}
\dot{V}&=e^T\dot{e}\\
&\le e^T[Ae+B\beta_1e+(D(y)\Theta-D(x)\hat{\Theta)}]\\
&\le e^T[Ae+B\beta_1e+(D(y)\hat{\Theta}-D(x)\hat{\Theta})-D(y)D^T(y)\lambda R^{-1}]\\
&\le e^T[Ae+B\beta_1e+\beta_2e\hat{\Theta}-D(y)D^T(y)W(y,e)R^{-1}]\\
&= e^TPe.
\end{split}
\end{equation}
From the TPBVP and Eq.\eqref{dl}, we know that $\lambda$ is the costate and can be described as a function of $y$ and $e$, here denoted by $W(y,e)$. By adjusting the  stable matrix $A_s$ and the function of horizon $T$, we aim to design a reasonable function $W(y, e)$ to make $\dot{V}$ smaller than zero
\begin{equation}
\dot{V}\le e^TPe \le 0,\quad \text{where} \quad P\le0.
\end{equation}
Note that $\dot{V}=0$ if and only if $e_i=0,$ $i=1,2,\cdots,n$. From the Barbalat's lemma \cite{astrom1995}, we can attain
\begin{equation}
 \lim_{t \rightarrow \infty}\|e_i\|=0, \, i=1,2,\cdots, n.
\end{equation}
It is clear that $e_i(t)\rightarrow 0$ as $t\rightarrow 0$. Thus, the synchronization error $e$ is asymptotically stable. Although the back-ward sweep algorithm is executable whenever the system is stable or not, with some choice of suitable stable matrix and horizon, we can also ensure the stability of the closed-loop nonlinear system by nonlinear receding horizon control.

\section{Application to chaotic systems}\label{sec:app2chaotic_sys}

In this section, we use a classical example from chaotic systems, namely the Lorenz oscillator \cite{lorenz1963}, to verify the effectiveness of the proposed method not only on time invariant parameters, but also on systems with time-varying parameters.

\subsection{Constant parameters}\label{sec:sec:app2chaotic_sys_const_params}

We first consider the parameter estimation problem of Lorenz chaotic system with constant parameters.

For this example, the Lorenz system is given by
\begin{equation}\label{system}
\begin{split}
&\frac{{\rm d}x(t)}{{\rm d}t}=\begin{bmatrix}10(x_2-x_1) \\ 28x_1-x_1x_3-x_2\\x_1x_2-\frac{8}{3}x_3\end{bmatrix},\\
&\frac{{\rm d}y(t)}{{\rm d}t}=\begin{bmatrix}\theta_1(x_2-x_1) \\ 28x_1-x_1x_3-x_2\\x_1x_2-\theta_2x_3\end{bmatrix}.
\end{split}
\end{equation}
where, the performance index is chosen as follows:
\begin{equation}
J=\int_t^{t+T}\{[y(\tau)-x(\tau)]^TQ[y(\tau)-x(\tau)]+[\Theta(\tau)-\hat{\Theta}(\tau)]^TR[\Theta(\tau)-\hat{\Theta}(\tau)]\}{\rm d}\tau.
\end{equation}
In this example, the weighting matrix $Q=R=\text{diag}(0.5,0.5,0.5)$.

The horizon $T$ in the performance index is given by
\begin{equation}
T(t)=T_f(1-e^{-\alpha t}),
\end{equation}
where $T_f=0.5$ and $\alpha=0.1$. It is clear that $T(t)$ satisfies $T(0)=0$ and $T(t)$ converges to the desired terminal time $T_f$ as time $t$ increases.

The stable matrix $A_s$ is chosen as $160I$. The initial states of the system are given by
\begin{equation}\label{ic_1}
\begin{bmatrix}x_1(0) \\ x_2(0)\\x_3(0)\end{bmatrix}=\begin{bmatrix}-3 \\ -3\\15\end{bmatrix},
\begin{bmatrix}y_1(0) \\ y_2(0)\\y_3(0)\end{bmatrix}=\begin{bmatrix}-10 \\ -10\\22\end{bmatrix}.
\end{equation}

\begin{figure}[htbp!]
\centering
   \includegraphics[width=12cm]{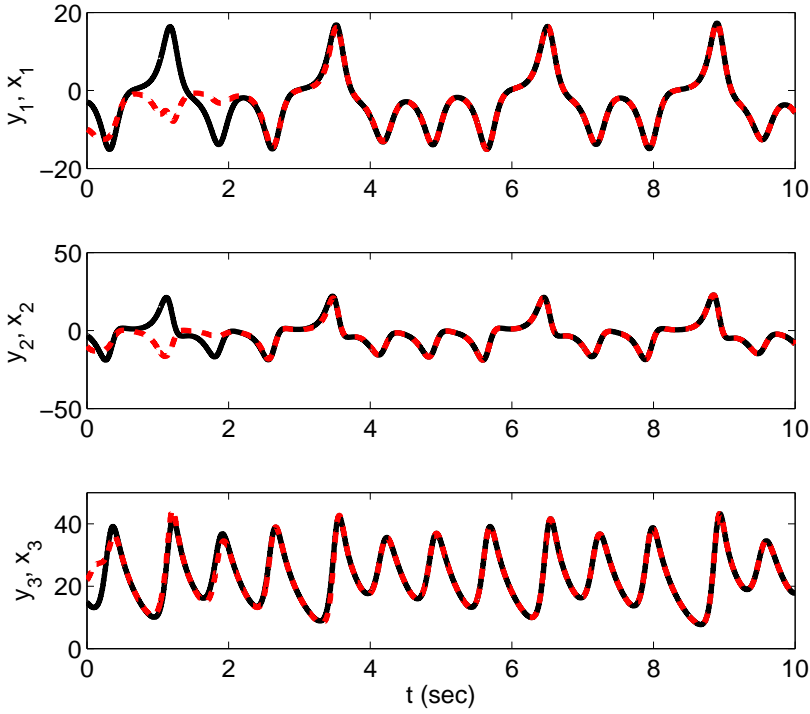}
   \vspace*{-.5cm}
   \caption{  (Color online) The trajectories of states with time invariant parameters where dash line denotes the trajectory of response system and solid line denotes the trajectory of reference system.}\label{fig1}
\end{figure}

\begin{figure}[htbp!]
\centering
   \includegraphics[width=12cm]{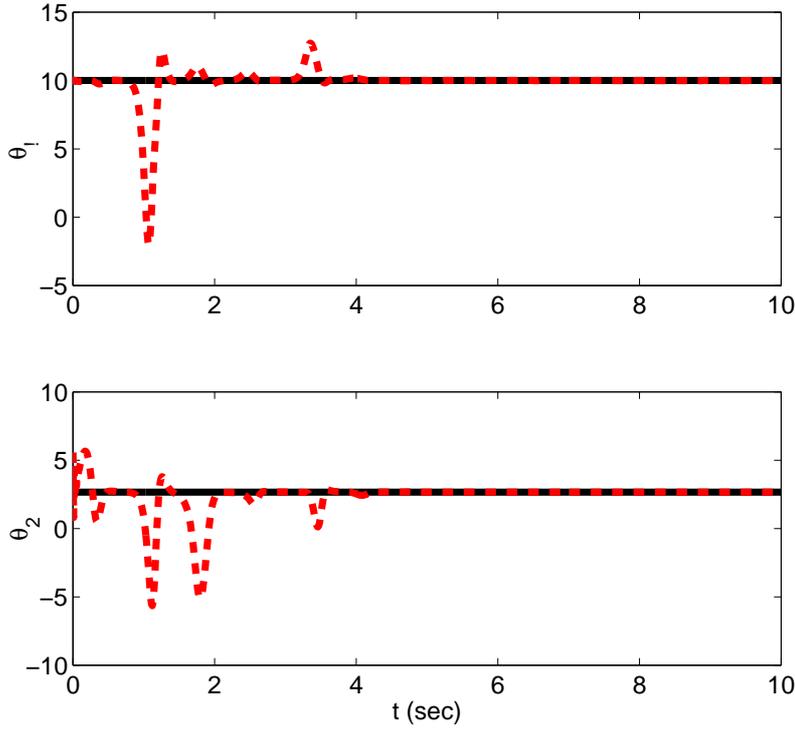}
   \vspace*{-.5cm}
   \caption{  (Color online) The trajectories of estimated constant parameters where dash line denotes the trajectory of estimated parameters and solid line denotes the trajectory of their true values.}\label{fig2}
\end{figure}


The simulation is implemented in MATLAB. The time step on the $t$ axis is $0.01$s and the time step on the $\tau$ axis is $0.005$s. Fig.\ref{fig1} shows the trajectories of drive-response systems in Eqs.\eqref{system} with initial conditions in Eqs.\eqref{ic_1}. The estimated parameters $\theta_1$ and $\theta_2$ are presented in Fig.\ref{fig2} which clearly show the estimated parameters converge to their true values by using the NRHC method.
%

\subsection{Time-varying parameters}\label{sec:sec:app2chaotic_sys_TimeVar_params}

In the following, we extend the NRHC estimation methodology to the case of time varying parameters. We consider the same Lorenz system, which was given in Section-\ref{sec:sec:app2chaotic_sys_const_params}, with time varying parameters $\theta_1=\frac{10\sin(t)}{t+1}$ and $\theta_2=\frac{8}{3}$.

The initial states are given by
\begin{equation}
\begin{bmatrix}x_1(0) \\ x_2(0)\\x_3(0)\end{bmatrix}=\begin{bmatrix}-3 \\ -3\\15\end{bmatrix},
\begin{bmatrix}y_1(0) \\ y_2(0)\\y_3(0)\end{bmatrix}=\begin{bmatrix}-6 \\ -6\\22\end{bmatrix}.
\end{equation}

\begin{figure}
\centering
   \includegraphics[width=12cm]{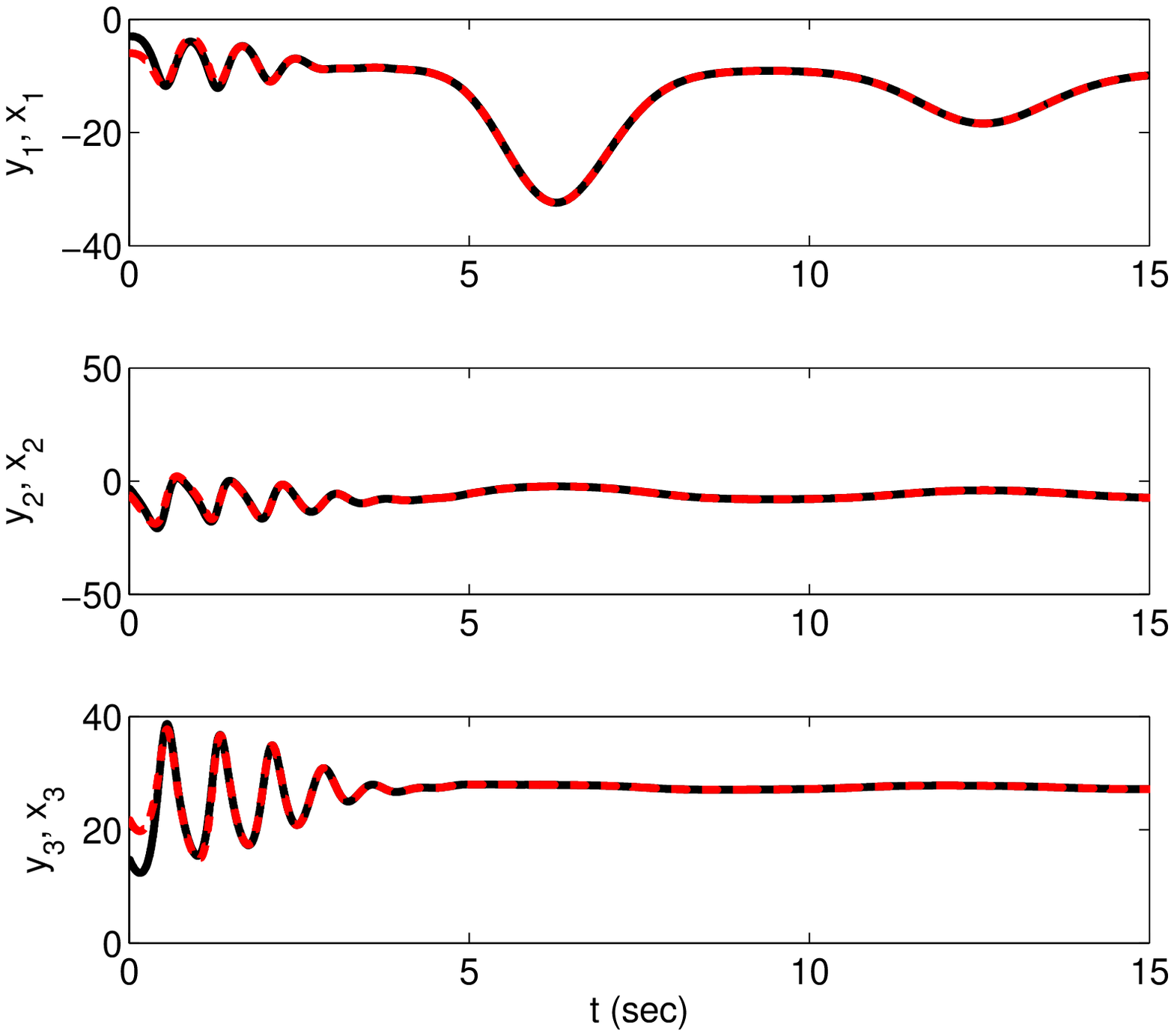}
   \caption{  (Color online) The trajectories of states with time varying parameters where dash line denotes the trajectory of response system and solid line denotes the trajectory of reference system.}\label{fig3}
   \includegraphics[width=12cm]{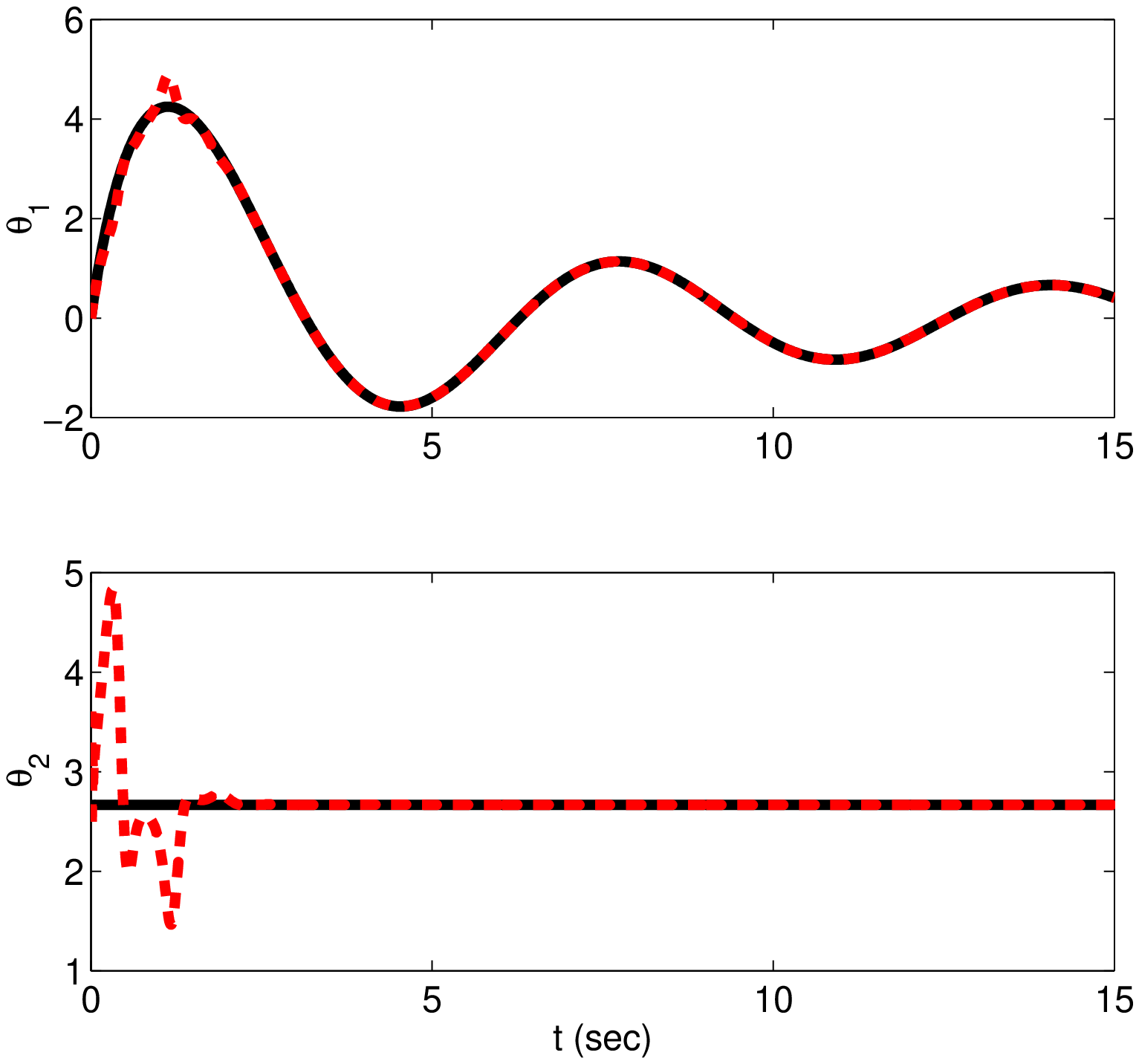}
   \caption{ (Color online) The trajectories of estimated time-varying parameters where dash line denotes the trajectory of estimated parameters and solid line denotes the trajectory of their true values.}\label{fig4}
\end{figure}

Fig.\ref{fig3} shows the trajectories of systems given in Eqs.\eqref{system} with time-varying parameters. The estimated parameters $\theta_1$ and $\theta_2$ are shown in Fig.\ref{fig4}. It is clear from Fig.\ref{fig3} and Fig.\ref{fig4} that in case of time varying characteristics NRHC still is able to perform as desired, and converges to the true values.

\subsection{Robustness against noise}\label{sec:sec:app2chaotic_sys_noise}

Noise, or generally speaking external disturbances, usually have significant effects on the performance and the outcomes of parameter estimation routine. Such external effects not only causes a drift in estimated parameters around the nominal value, but also results in potentially unstable systems. In the following, we investigate the effect of the noise in aforementioned dynamics.

\begin{figure}[htbp!]
 \includegraphics[width=11cm]{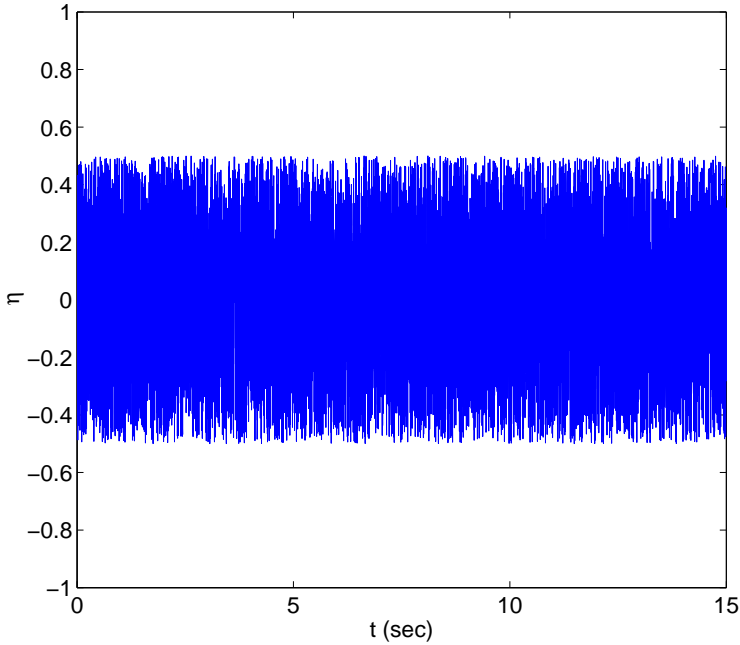}
   \caption{ (Color online) The trajectories of noise signal $\eta(t)$ in the time invariant systems.}\label{fig5}
\end{figure}

\begin{figure}[htbp!]
\centering
   \includegraphics[width=12cm]{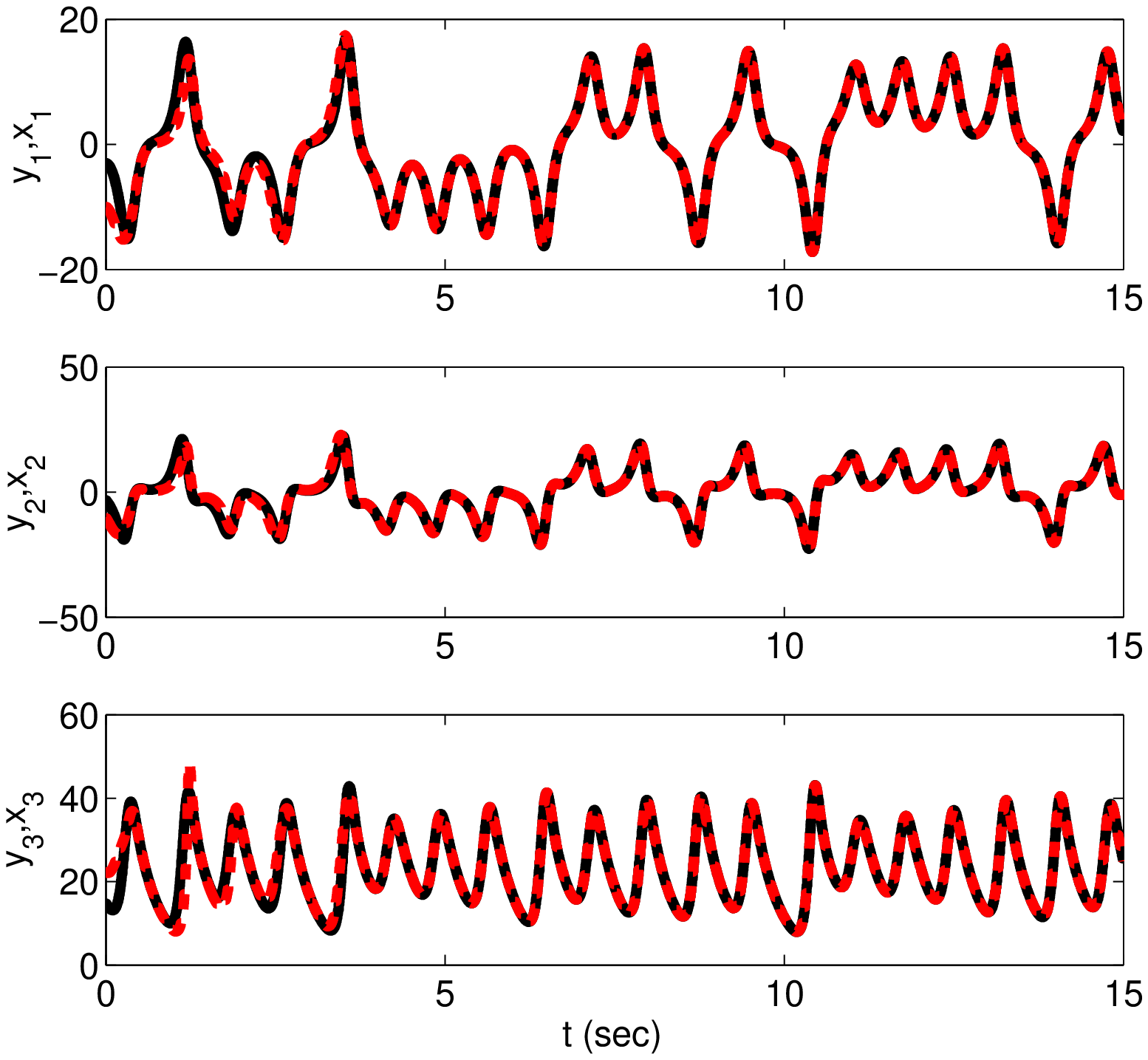}
   \caption{  (Color online) The trajectories of states with constant parameters in presence of noise where dash line denotes the trajectory of response system and solid line denotes the trajectory of reference system.}\label{fig6}
\end{figure}

\begin{figure}[htbp!]
\centering
   \includegraphics[width=12cm]{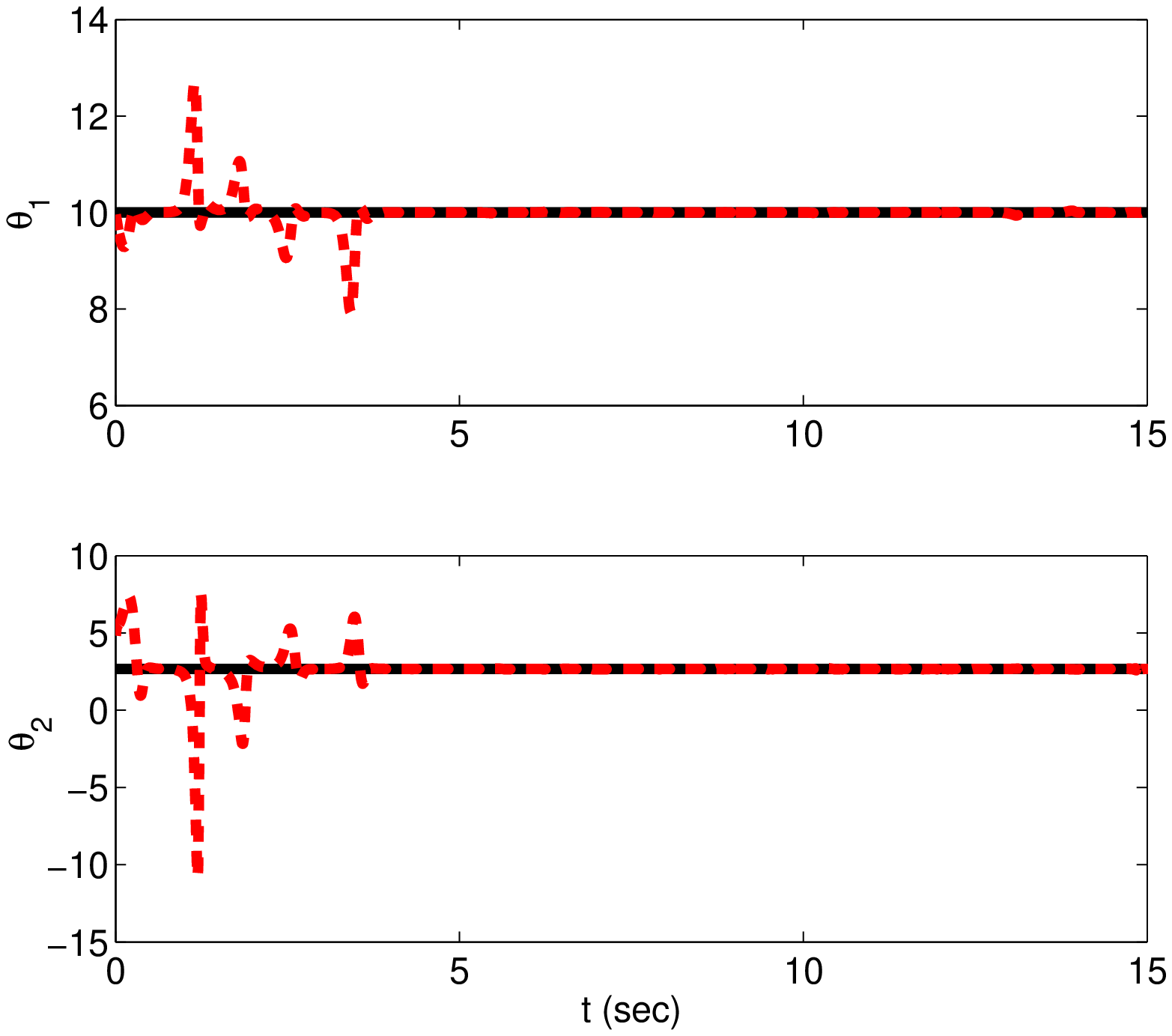}
   \caption{ (Color online) The trajectories of estimated constant parameters in presence of noise where dash line denotes the trajectory of estimated parameters and solid line denotes the trajectory of their true values.}\label{fig7}
\end{figure}

We first consider the case where the noise propagates in the drive system with constant parameters which could be expressed as
\begin{equation}\label{system_n}
\begin{split}
&\frac{{\rm d}x(t)}{{\rm d}t}=\begin{bmatrix}10(x_2-x_1)+\eta(t) \\ 28x_1-x_1x_3-x_2+\eta(t)\\x_1x_2-\frac{8}{3}x_3+\eta(t)\end{bmatrix},\\
&\frac{{\rm d}y(t)}{{\rm d}t}=\begin{bmatrix}\theta_1(x_2-x_1) \\ 28x_1-x_1x_3-x_2\\x_1x_2-\theta_2x_3\end{bmatrix},
\end{split}
\end{equation}
where $\eta(t)$ represents the band-limited white noise. The simulation results are depicted in Figs.\ref{fig5},\ref{fig6},\ref{fig7}. In
this case, the estimated constant parameters precisely match their original values, which demonstrates the robustness characteristic of proposed, NRHC based, parameter estimation routine.


\begin{figure}
 \includegraphics[width=12cm]{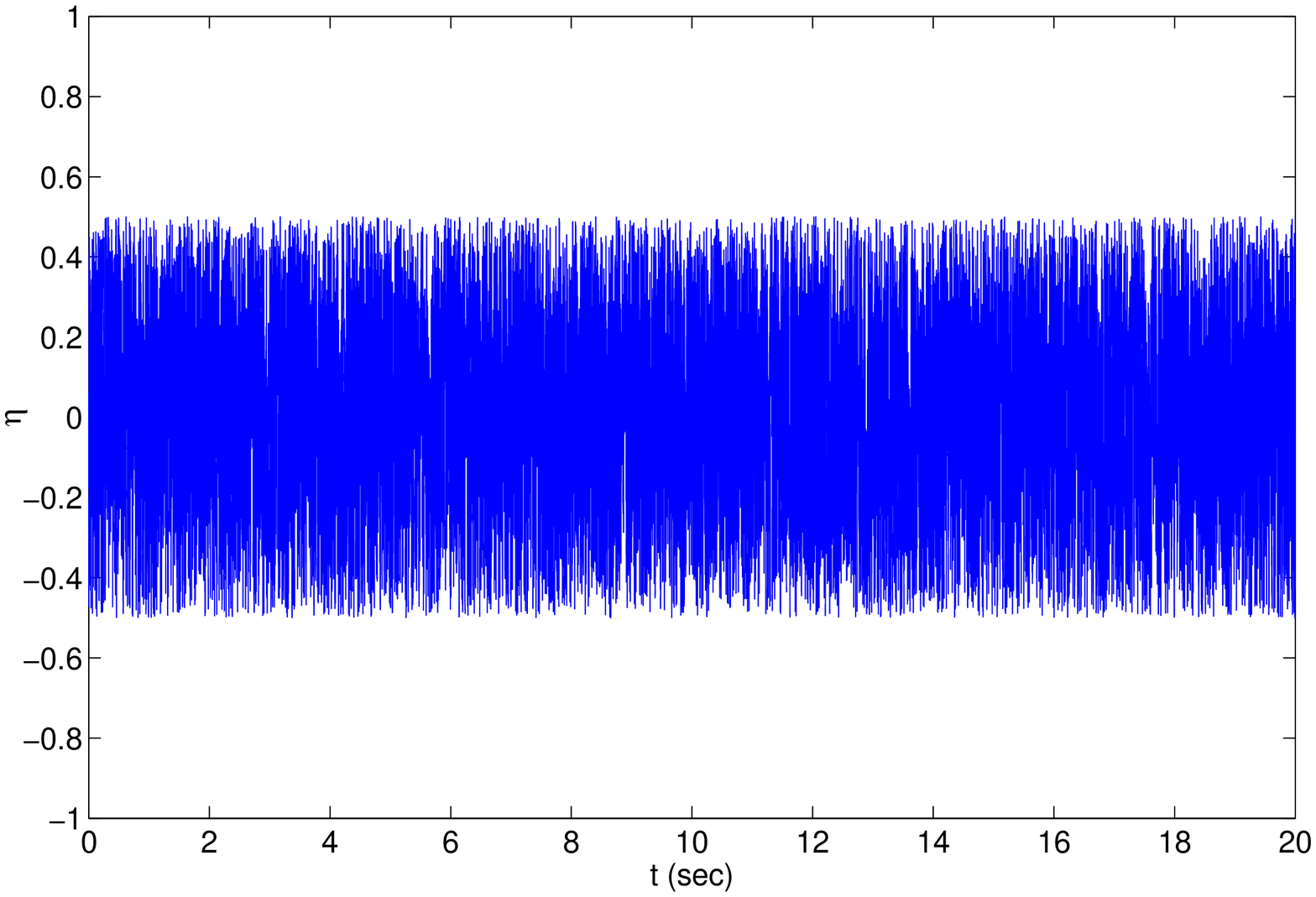}
   \caption{ (Color online) The trajectories of noise signal $\eta(t)$ in the time-varying systems.}\label{fig8}
\end{figure}

\begin{figure}[!htb]
\centering
   \includegraphics[width=12cm]{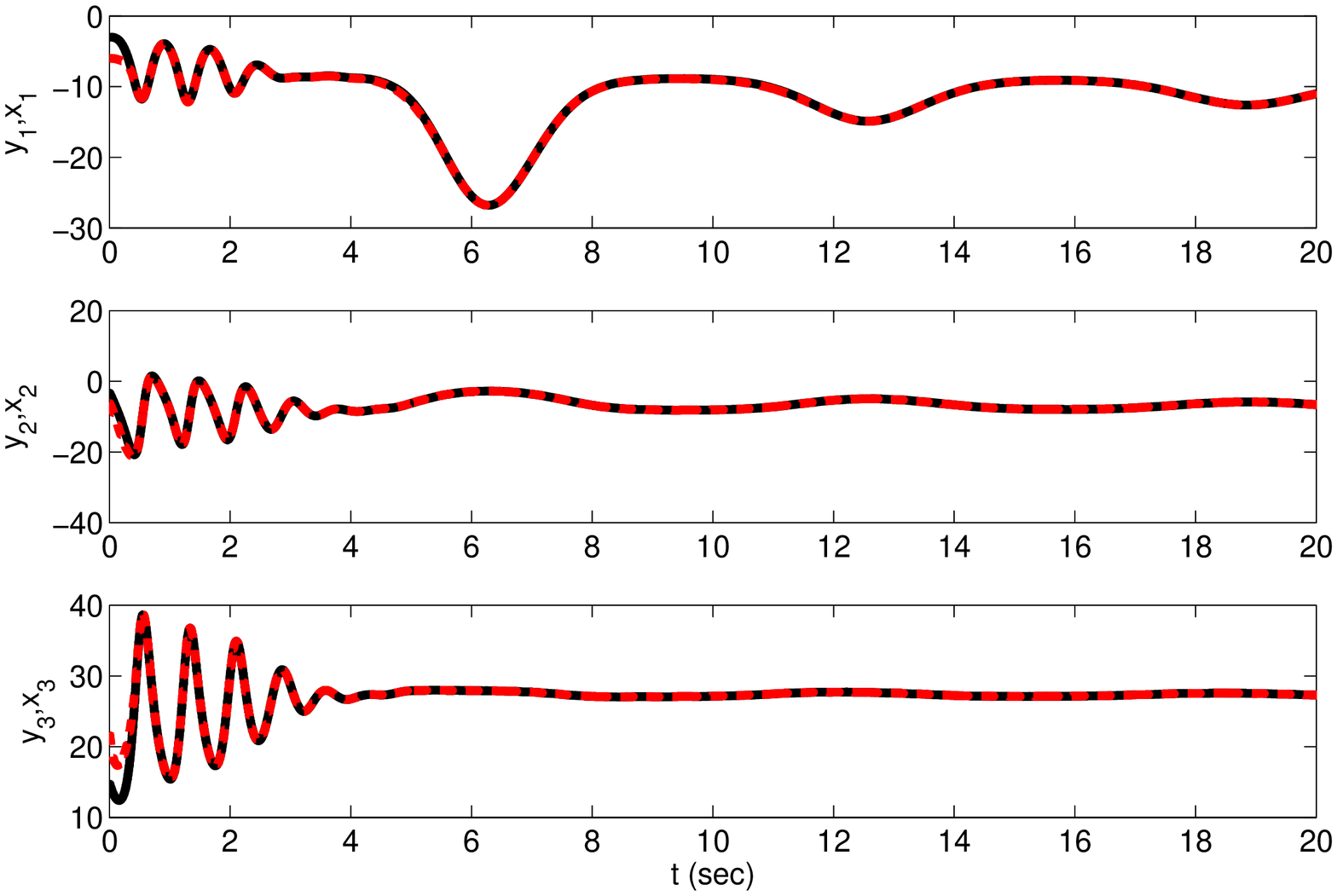}
   \caption{  (Color online) The trajectories of states with time varying parameters in presence of noise where dash line denotes the trajectory of response system and solid line denotes the trajectory of reference system.}\label{fig9}
\end{figure}

\begin{figure}[!htb]
\centering
   \includegraphics[width=12cm]{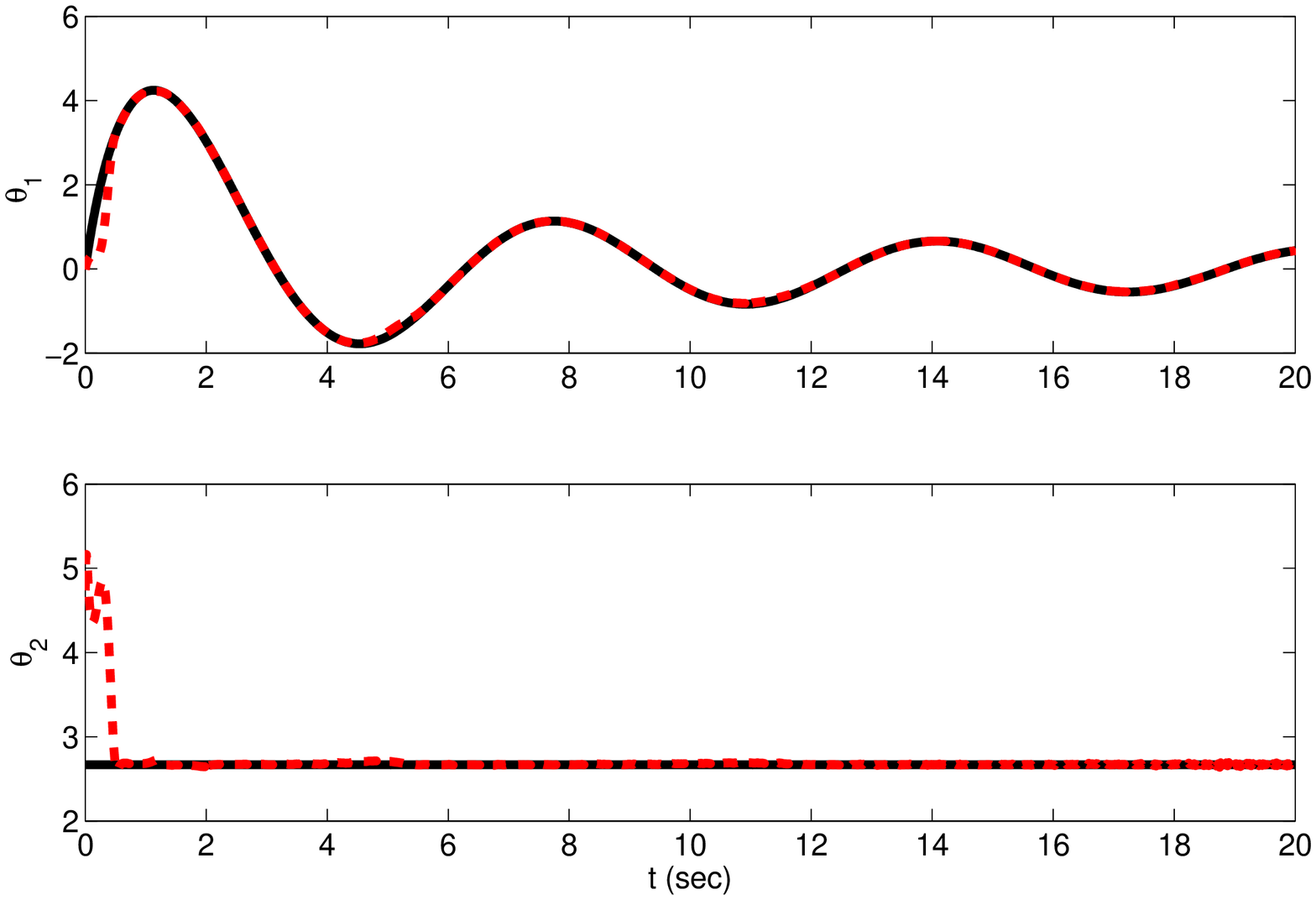}
   \caption{ (Color online) The trajectories of estimated time-varying parameters in presence of noise where dash line denotes the trajectory of estimated parameters and solid line denotes the trajectory of their true values.}\label{fig10}
\end{figure}

Next, we propagate the noise content in the drive system with time-varying parameters which
is expressed as

\begin{equation}\label{system_n}
\begin{split}
&\frac{{\rm d}x(t)}{{\rm d}t}=\begin{bmatrix}\frac{10\sin(t)}{t+1}(x_2-x_1)+\eta(t) \\ 28x_1-x_1x_3-x_2+\eta(t)\\x_1x_2-\frac{8}{3}x_3+\eta(t)\end{bmatrix},\\
&\frac{{\rm d}y(t)}{{\rm d}t}=\begin{bmatrix}\theta_1(x_2-x_1) \\ 28x_1-x_1x_3-x_2\\x_1x_2-\theta_2x_3\end{bmatrix},
\end{split}
\end{equation}
where again $\eta(t)$ is the band-limited white noise. The simulation results are illustrated in Figs.\ref{fig8},\ref{fig9},\ref{fig10}. In
this case, the estimated time-varying parameters also demonstrate a good match with their original values, in presence of noise.



\section{Discussions \& Conclusions}\label{sec:conclusion}
In this paper, a novel method based on real time nonlinear receding horizon control is proposed for estimating unknown parameters of general nonlinear and chaotic systems. In this specific set-up, the estimation problem is reduced to a form of solving the nonlinear receding horizon optimization problem as a parameter optimization method. Based on the stabilized continuation method, the back-ward sweep algorithm is introduced to integrate the costate in real time and to minimize the estimation error. The algorithm does not require any stability assumption of the system and also can guarantee the stability with some suitable choice of stable matrix and horizon length. The method is applicable for both time invariant and time varying dynamics with noise, which demonstrates the power of the methodology.





\end{document}